\documentclass[11pt,a4paper]{article}
\usepackage[T1]{fontenc}
\usepackage[utf8]{inputenc}
\usepackage{lmodern}
\usepackage[margin=27mm]{geometry}
\usepackage{amsmath,amssymb,amsthm}
\usepackage{microtype}
\usepackage{enumitem}
\usepackage[hidelinks]{hyperref}
\hypersetup{pdftitle={Sharp bounds for Duke--Erd\H{o}s--R\"odl problems on cycle-connected subgraphs},pdfauthor={Guorong Gao and Run Zhao},pdfsubject={The Duke--Erdos--Rodl problem on cycle-connected subgraphs},pdfkeywords={cycle-connected subgraph, short even cycle, extremal graph theory, biaffine incidence graph}}
\newtheorem{theorem}{Theorem}[section]
\newtheorem{lemma}[theorem]{Lemma}
\newtheorem{problem}[theorem]{Problem}
\newtheorem{proposition}[theorem]{Proposition}

\setlist[enumerate]{leftmargin=*,itemsep=2pt,topsep=4pt}
\title{\bfseries Sharp bounds for Duke--Erd\H{o}s--R\"odl problems on cycle-connected subgraphs}
\author{Guorong Gao\thanks{Corresponding author. E-mail address: \texttt{grgao@fzu.edu.cn}.}\quad and\quad Run Zhao\thanks{E-mail address: \texttt{runzhao0628@163.com}.}\\[5pt]
\small School of Mathematics and Statistics, Fuzhou University,\\
\small Fuzhou, Fujian, China}
\date{}
\begin{document}
\maketitle

\begin{abstract}
For $k\ge2$, a graph is $C_{2k}$-connected if every two distinct edges lie together on an even cycle of length at most $2k$ in that graph. Let $f_k(n,m)$ be the largest number of edges in a $C_{2k}$-connected subgraph that every $n$-vertex graph with $m$ edges must contain. In this paper, we study two cycle-connected problems raised by Duke, Erd\H{o}s and R\"odl. 
\begin{itemize}
    \item  First, Duke, Erd\H{o}s and R\"odl showed that $f_3(n,\lfloor n^{2-\beta}\rfloor)=\Theta(n^{2-3\beta})$ for every fixed $0< \beta<1/2$ and ask for the case of $\beta\ge 1/2$. We solve this problem in a stronger form by determining that the optimal range of $\beta$ is $(0,2/3]$. 
    \item Second, we prove that $f_2(n,\binom n2-s)=(1/(2e)+o(1))n^3/s$ whenever $\frac{s}{n}\to\infty$ and $s=o(n^2)$, determining the leading constant in an estimate of Duke, Erd\H{o}s, and R\"odl. 
\end{itemize}
The first result uses a finite-field construction and a count of three-edge paths. The second uses the known structure of $C_4$-connected graphs and estimates for complete bipartite subgraphs in a random graph.
\end{abstract}

\section{Introduction}\label{sec:introduction}
All graphs are finite and simple, and subgraphs need not be induced. We write $v(G)$ and $e(G)$ for the numbers of vertices and edges of $G$, $N_G(x)$ for the neighborhood of $x$, $d_G(x)=|N_G(x)|$, and $\delta(G)$ for the minimum degree. The diameter of a disconnected graph is infinite. As usual, $K_n$, $K_{a,b}$, and $C_\ell$ denote the complete graph, the complete bipartite graph with part sizes $a,b$, and the cycle of length $\ell$.

For an integer $k\ge2$, a graph $H$ is \emph{$C_{2k}$-connected} if every two distinct edges lie together on a cycle in $H$ belonging to $\{C_4,C_6,\ldots,C_{2k}\}$. Thus a $C_6$-connected graph may use either a four-cycle or a six-cycle. A one-edge graph satisfies the definition. For $0\le m\le\binom n2$, put
\begin{equation}\label{eq:minimax}
f_k(n,m):=\min_{\substack{v(G)=n\\ e(G)=m}}\max\{e(H):H\subseteq G\text{ is }C_{2k}\text{-connected}\}.
\end{equation}
Thus $f_k(n,m)$ is the largest number of edges in a $C_{2k}$-connected subgraph that every $n$-vertex graph with $m$ edges must contain.
We write $[t]=\{1,\ldots,t\}$ and use $\mathbb F_q$ for the field of order $q$. All logarithms are natural, and asymptotic notation refers to $n\to\infty$. Constants with a subscript $\beta$ may depend on the fixed exponent $\beta$.

Duke and Erd\H{o}s~\cite{DukeErdos} introduced the problem of finding large subgraphs whose pairs of edges lie on short common cycles. Duke, Erd\H{o}s, and R\"odl developed this problem further in~\cite{DER1984,DER1992}. Their results give $f_3(n,\lfloor n^{2-\beta}\rfloor)=\Theta_\beta(n^{2-3\beta})$ for $0< \beta< 1/2$. Fox and Sudakov~\cite{FoxSudakov} later showed that $f_4(n,\lfloor n^{2-\beta}\rfloor)=\Theta_\beta(n^{2-2\beta})$ for $0< \beta< 1/5$. In a recent preprint, Li~\cite{Li} improved the range of $\beta$ to $0< \beta\le 1/3$. A natural problem proposed by Duke, Erd\H{o}s, and R\"odl \cite{DER1992} asks for the behavior of $f_3(n,\lfloor n^{2-\beta}\rfloor)$ for $1/2\leq \beta\leq 1 $.

\begin{problem}
Is $f_3(n,\lfloor n^{2-\beta}\rfloor)$ unbounded for $1/2\leq \beta\leq 1 $?
\end{problem}

Our first result solves this problem in a stronger form.
For completeness, the endpoint $\beta=0$ is handled separately: the expression $\lfloor n^2\rfloor$ lies outside the domain $m\le\binom n2$, while the result of ~\cite{DukeErdos} implies that $f_3(n,cn^2)=\Theta(n^2)$ for every constant $0<c<1/2$. Theorem~\ref{thm:sparse} therefore states the nontrivial power-law range $0<\beta\le2/3$.

\begin{theorem}\label{thm:sparse}
For every fixed real number $\beta$ with $0<\beta\le2/3$,
\[
f_3(n,\lfloor n^{2-\beta}\rfloor)=\Theta_\beta(n^{2-3\beta}).
\]
\end{theorem}

At $\beta=2/3$, our proof gives $1\le f_3(n,\lfloor n^{4/3}\rfloor)\le580800$ for all sufficiently large $n$. For fixed $2/3<\beta\leq 1$, the classical girth-eight constructions of Benson~\cite{Benson} give $f_3(n,\lfloor n^{2-\beta}\rfloor)=1$ for large $n$. Indeed, these constructions have $\Theta(q^3)$ vertices and $\Theta(q^4)$ edges; choosing a suitable prime $q=\Theta(n^{1/3})$, deleting edges, and adding isolated vertices gives a host graph with no four-cycle or six-cycle. Thus, Theorem~\ref{thm:sparse} determines the order for every fixed $0\le\beta \le1$.

Our second result concerns $C_4$-connected subgraphs in near-complete graphs. Duke, Erd\H{o}s, and R\"odl~\cite[Theorem~2]{DER1992} proved that, whenever $s/n\to\infty$ and $s=o(n^2)$,
\[
\left(\frac1{16}-o(1)\right)\frac{n^3}{s}
\le f_2\left(n,\binom n2-s\right)
\le (1+o(1))\frac{n^3}{s}.
\]
They left the asymptotic value open~\cite[p.~268]{DER1992}. We determine its leading constant throughout the same range.

\begin{theorem}\label{thm:dense}
Let $s=s(n)$ be an integer satisfying $s/n\to\infty$ and $s=o(n^2)$. Then
\[
f_2\left(n,\binom n2-s\right)=\left(\frac{1}{2e}+o(1)\right)\frac{n^3}{s}.
\]
Moreover, every host graph contains a complete bipartite subgraph that gives this asymptotic lower bound.
\end{theorem}

The complete multipartite structure of $C_4$-connected graphs and the fixed-density estimates were established in~\cite[Proposition~0 and Theorem~1]{DER1992}; related complete bipartite subgraph problems appear in~\cite{EFRS}. We use the counting and random-graph methods of~\cite{DER1992}, with estimates that remain valid as the density tends to one. This yields the constant in Theorem~\ref{thm:dense}.

Section~\ref{sec:structure} proves Theorem~\ref{thm:sparse}. Its upper bound uses a finite-field incidence graph with a product-form edge equation, and its lower bound counts three-edge paths with common endpoints. Section~\ref{sec:dense} then proves Theorem~\ref{thm:dense}.

\section{Proof of Theorem~\ref{thm:sparse}}\label{sec:structure}
\subsection{Bipartite structure}
We first describe bipartite $C_6$-connected graphs. This description will be used in the construction and the lower-bound proof.

\begin{lemma}\label{lem:bipartite}
Let $H=(X,Y;E)$ be a bipartite graph with no isolated vertices and at least two edges. Then $H$ is $C_6$-connected if and only if
\[
\delta(H)\ge 2\qquad\text{and}\qquad \operatorname{diam}(H)\le 3.
\]
Equivalently, $\delta(H)\ge 2$ and every two distinct vertices in the same bipartition class have a common neighbor.
\end{lemma}

\begin{proof}
Assume first that $H$ is $C_6$-connected. Every edge lies on a common cycle with another edge, so $\delta(H)\ge2$.

For distinct vertices $u,v$ in the same class, choose an edge incident with each. These edges lie on a common $C_4$ or $C_6$, in which $u,v$ have a common neighbor. Hence every such pair has distance two.

If $x\in X$ and $y\in Y$ are nonadjacent, choose $x'\in N_H(y)$; a common neighbor of $x,x'$ gives an $x$--$y$ path of length three. Thus the diameter is at most three. Conversely, this diameter bound forces distinct vertices in the same class to have distance two, proving the equivalence of the two formulations.

Conversely, suppose that $\delta(H)\ge 2$ and that any two distinct vertices in the same class have a common neighbor. Let $e=xy$ and $f=x'y'$ be distinct edges, where $x,x'\in X$ and $y,y'\in Y$.

First suppose that the edges share an endpoint, say $x=x'$. Choose $x_1\in N_H(y)\setminus\{x\}$ and $x_2\in N_H(y')\setminus\{x\}$. If $x_1=x_2$, then $x-y-x_1-y'-x$ is a $C_4$ containing both edges. Otherwise, choose $z\in N_H(x_1)\cap N_H(x_2)$. If $z\notin\{y,y'\}$, then $x-y-x_1-z-x_2-y'-x$ is a $C_6$. If $z=y$, then $x-y-x_2-y'-x$ is a $C_4$, while if $z=y'$, then $x-y-x_1-y'-x$ is a $C_4$. The case in which the common endpoint belongs to $Y$ is symmetric.

We may therefore suppose that $x\ne x'$ and $y\ne y'$. If both cross-edges $xy'$ and $x'y$ are present, then $x-y-x'-y'-x$ is a $C_4$. If neither cross-edge is present, choose $a\in N_H(x)\cap N_H(x')$ and $b\in N_H(y)\cap N_H(y')$. The missing cross-edges imply $a\notin\{y,y'\}$ and $b\notin\{x,x'\}$, so $x-y-b-y'-x'-a-x$ is a simple $C_6$ containing $e$ and $f$.

It remains to consider the case in which exactly one cross-edge is present. By symmetry, assume that $xy'\in E(H)$ but $x'y\notin E(H)$. Choose $a\in N_H(x')\setminus\{y'\}$ and $b\in N_H(y)\setminus\{x\}$. Notice that $a\ne y$ and $b\ne x'$ because $x'y$ is absent. If $xa\notin E(H)$, choose $z\in N_H(a)\cap N_H(y)$. Then $z\notin\{x,x'\}$, so $x-y-z-a-x'-y'-x$ is a simple $C_6$. If $xa\in E(H)$ but $by'\notin E(H)$, choose $w\in N_H(b)\cap N_H(x')$. Now $w\notin\{y,y'\}$, and $x-y-b-w-x'-y'-x$ is a simple $C_6$. Finally, if both $xa$ and $by'$ are present, then $x-y-b-y'-x'-a-x$ is a simple $C_6$. Thus any two distinct edges lie on a common $C_4$ or $C_6$.
\end{proof}

Let $B_Q$ be the biaffine incidence graph over $K=\mathbb{F}_Q$. Its two vertex classes represent points and nonvertical lines. The point vertices are $(x,y)\in K^2$, and the line vertices are $[s,a]\in K^2$. A point and a line are adjacent exactly when the point lies on the line:
\[
(x,y)\sim[s,a]\quad\Longleftrightarrow\quad y=sx+a.
\]
The standard incidence properties of the affine plane~\cite{Hirschfeld} show that $B_Q$ is $Q$-regular, has $2Q^2$ vertices and $Q^3$ edges, and is $C_4$-free: each point lies on one line of each slope, and two distinct points lie on at most one common nonvertical line.

We call a bipartite graph \emph{pair-covered} if every pair of distinct vertices in the same bipartition class has a common neighbor. We will use finite point--line incidence structures with the following properties: every two distinct points lie on exactly one line, every line contains at least two points, and there are at least two lines. Here a line is simply a subset of the point set. A \emph{near-pencil} is such a structure with one line containing all but one point. The remaining point is called the exceptional point. Each other line joins the exceptional point to one point on the long line and contains exactly these two points.

\begin{lemma}\label{lem:support}
Let $Q=r^2$, where $r$ is prime. Suppose that $J\subseteq B_Q$ has at least one edge and no isolated vertices, and that every pair of distinct vertices in either bipartition class has a common neighbor in $J$. Then $e(J)\le 3Q$.

Moreover, if $\delta(J)\ge 2$, then $J$ is the incidence graph of a near-pencil. More explicitly, for some point--line pair $(P_0,\ell_0)$ with the point not on the line, it is the union of $u$ three-edge paths from $P_0$ to $\ell_0$ that share no internal vertices, and $e(J)=3u$.
\end{lemma}

\begin{proof}
Let $\mathcal{P}$ and $\mathcal{L}$ be the point and line vertices of $J$. We call them the support of $J$.

\smallskip\noindent\emph{All incidences between the selected points and lines.}
We first show that $J$ contains every incidence of $B_Q$ between $\mathcal{P}$ and $\mathcal{L}$. Suppose that $P\in\mathcal{P}$ and $\ell\in\mathcal{L}$ are incident in $B_Q$ but $P\ell\notin E(J)$. Choose $P'\in N_J(\ell)$. Then $P'\ne P$. By pair-covering, there is a line $\ell'\in\mathcal{L}$ incident in $J$ with both $P$ and $P'$. Necessarily $\ell'\ne\ell$, and $P,\ell,P',\ell'$ form a $C_4$ in $B_Q$, a contradiction. Hence $J$ is the incidence graph induced by $\mathcal{P}\cup\mathcal{L}$ inside $B_Q$.

\smallskip\noindent\emph{The degree-one case.}
Suppose that a point $P$ has the unique neighbor $\ell$. Pair-covering forces every selected point to lie on $\ell$. Each other selected line has exactly one selected point, by $C_4$-freeness and the absence of isolated vertices. If such a line exists, let its selected point be $P_0$. Pair-covering on the line side forces every line other than $\ell$ to contain $P_0$. There are at most $Q-1$ such lines and at most $Q$ points on $\ell$, so $e(J)\le2Q-1$. If $\ell$ is the only line, then $e(J)\le Q$. Exchanging points and lines gives the same bound when a selected line has degree one.

\smallskip\noindent\emph{All degrees are at least two.}
It remains to assume that $\delta(J)\ge 2$. Pair-covering and $C_4$-freeness imply that any two selected points lie on exactly one selected line and any two selected lines meet in exactly one selected point. Also, each selected line contains at least two selected points, and each selected point lies on at least two selected lines. Thus the point--line incidence structure and its dual, obtained by exchanging points and lines, both satisfy the properties stated above. In particular, both have at least two lines. No selected line contains all selected points: a second selected line would contain at least two of those points, contrary to $C_4$-freeness. The same argument applies to the dual.

Put $v=|\mathcal{P}|$ and $b=|\mathcal{L}|$. The de Bruijn--Erd\H{o}s theorem gives $b\ge v$ for a finite point--line incidence structure with these properties. Its equality cases are precisely near-pencils and finite projective planes~\cite[Theorem~1]{deBruijnErdos}. Applying the same theorem to the dual gives $v\ge b$. Hence $v=b$, and the structure is a near-pencil or a finite projective plane.

We exclude the projective-plane case. Embed the biaffine incidence structure in $\operatorname{PG}(2,K)$ by adding the vertical affine lines, all points at infinity, and the line at infinity. The line joining two selected points is their unique selected line, and the intersection of two selected lines is their unique selected point. We have also shown that the selected points and lines inherit every incidence from the ambient plane. Thus, in the projective-plane case, $(\mathcal{P},\mathcal{L})$ is a projective subplane of $\operatorname{PG}(2,r^2)$.

The subplanes of $\operatorname{PG}(2,p^h)$ have order $p^d$, where $d$ is a positive divisor of $h$~\cite[Remark~14.3.51]{EbertStorme}. Since $r$ is prime and $h=2$, the order $q_0$ of this subplane is either $r$ or $r^2$.

On the other hand, no two points of $\mathcal{P}$ can have the same first coordinate: the line joining them would be vertical, while $\mathcal{L}$ contains only nonvertical lines. Consequently, $|\mathcal{P}|\le Q=r^2$. But a projective plane of order $q_0$ has $q_0^2+q_0+1$ points, and $q_0^2+q_0+1\ge r^2+r+1>Q$, a contradiction.

Thus the incidence structure is a near-pencil. Let $P_0$ be its exceptional point, let $\ell_0$ be its long line, and let $P_1,\ldots,P_u$ be the ordinary points on $\ell_0$. Its remaining lines are precisely the two-point lines $\ell_i=P_0P_i$, $1\le i\le u$. Consequently, its incidence graph is the union, for $1\le i\le u$, of the internally vertex-disjoint paths $P_0-\ell_i-P_i-\ell_0$. The long line contains at most $Q$ points in $B_Q$, so $u\le Q$ and $e(J)=u+2u=3u\le 3Q$.
\end{proof}

\subsection{The finite-field construction}\label{sec:construction}
We now construct a graph using copies of points and lines in $B_Q$. Lemma~\ref{lem:support} describes the possible supports of its $C_6$-connected subgraphs, and the next intersection bound controls their size.

Fix an odd prime $r>4$, put $Q=r^2$, and write $K=\mathbb{F}_{r^2}=\mathbb{F}_r(\omega)$ with $\omega^2=\eta$, where $\eta\in\mathbb{F}_r$ is a nonsquare. For $\alpha,\delta\in\mathbb{F}_r$, define
\[
L_\alpha=\{u+\omega\alpha:u\in\mathbb{F}_r\},
\qquad
P_\delta=\{v+\omega(v^2+\delta):v\in\mathbb{F}_r\}.
\]
Each of the two families partitions $K$. If $z=z_0+\omega z_1\in K$, then $z$ belongs to the unique horizontal layer $L_{z_1}$ and to the unique parabolic layer $P_{z_1-z_0^2}$.

We next define two periodic sequences. Choose an ordering $\zeta_1,\ldots,\zeta_r$ of $\mathbb{F}_r$. List the elements $\zeta_p+\omega\zeta_q$, where $q=1,\ldots,r$ and $p=1,\ldots,r$, in lexicographic order of $(q,p)$, with $p$ varying first. Repeat this length-$Q$ list to obtain the sequence $\lambda_1,\lambda_2,\ldots$. Similarly, list the elements $\zeta_p+\omega(\zeta_p^2+\zeta_q)$ in the same order and repeat the list to obtain the sequence $\xi_1,\xi_2,\ldots$. Every element of $K$ occurs exactly once in each period of each sequence. Thus, among the first $T$ terms of either sequence, every element occurs at most $\lceil T/Q\rceil$ times, and the terms are distinct when $T\le Q$. Moreover, the first $T$ terms of the $\lambda$-sequence meet at most $\lceil T/r\rceil$ horizontal layers, while those of the $\xi$-sequence meet at most $\lceil T/r\rceil$ parabolic layers.

For an integer $t\ge r$, define $G(r,t)$ with two disjoint vertex classes indexed by
\[
\{(i,a):i\in[t],\ a\in K\}
\quad\text{and}\quad
\{(j,b):j\in[t],\ b\in K\},
\]
where
\begin{equation}\label{eq:lift}
(i,a)(j,b)\in E(G(r,t))\quad\Longleftrightarrow\quad b=\lambda_i\xi_j+a.
\end{equation}
Each left vertex has one neighbor for every $j\in[t]$, and each right vertex has one neighbor for every $i\in[t]$. Thus $G(r,t)$ is $t$-regular, with
\begin{equation}\label{eq:counts}
v(G(r,t))=2tQ,\qquad e(G(r,t))=t^2Q.
\end{equation}
Geometrically, $(i,a)$ is a copy of the line $[\lambda_i,a]$, and $(j,b)$ is a copy of the point $(\xi_j,b)$ in $B_Q$. 

Write $K^\times=K\setminus\{0\}$ for the multiplicative group of $K$.

\begin{lemma}\label{lem:intersection}
Let $S,X\subseteq K$, where $S$ is contained in the union of $h$ horizontal layers $L_\alpha$ and $X$ is contained in the union of $k$ parabolic layers $P_\delta$. For all $s_0,x_0\in K$ and $c\in K^\times$,
\[
\bigl|\{(s,x)\in S\times X:(s-s_0)(x-x_0)=c\}\bigr|\le 4hk.
\]
\end{lemma}

\begin{proof}
There are at most $hk$ relevant layer pairs, so it is enough to bound the number of solutions in each $L_\alpha\times P_\delta$ by four.

Write $s_0=u_0+\omega\alpha_0$ and $x_0=v_0+\omega(v_0^2+\delta_0)$; these representations are unique. Let $s=u+\omega\alpha\in L_\alpha$ and $x=v+\omega(v^2+\delta)\in P_\delta$. Set $U=u-u_0$, $w=v-v_0$, $A=\alpha-\alpha_0$, $B=\delta-\delta_0$, and $\Psi(w)=w^2+2v_0w+B$. Then
\[
s-s_0=U+\omega A,\qquad x-x_0=w+\omega\Psi(w).
\]
Write $c=c_0+\omega c_1$. Then $(c_0,c_1)\ne(0,0)$ because $c\ne0$. The product equation is equivalent to
\[
Uw+\eta A\Psi(w)=c_0,\qquad U\Psi(w)+Aw=c_1.
\]
Multiply the first equation by $\Psi(w)$ and the second by $w$, and then subtract. Every solution therefore gives a root of
\begin{equation}\label{eq:polynomial}
\Phi(w):=\eta A\Psi(w)^2-Aw^2-c_0\Psi(w)+c_1w=0.
\end{equation}
The polynomial in~\eqref{eq:polynomial} has degree at most 4 and hence has at most four roots in $\mathbb{F}_r$.

Each root $w$ determines $x\in P_\delta$. A solution has $x\ne x_0$, and then $s=s_0+c(x-x_0)^{-1}$ is unique. Hence each layer pair contributes at most four solutions, proving the claim.
\end{proof}

\paragraph{The construction bound.}
The next proposition bounds the number of edges in each $C_6$-connected subgraph of $G(r,t)$.

\begin{proposition}\label{prop:construction}
Let $r>4$ be an odd prime, let $Q=r^2$, and let $t\ge r$ be an integer. There is an explicit bipartite graph $G(r,t)$ satisfying
\[
v(G(r,t))=2tQ,\qquad e(G(r,t))=t^2Q,
\]
such that every $C_6$-connected subgraph $H\subseteq G(r,t)$ satisfies $e(H)\le48t^2/Q$.
\end{proposition}

\begin{proof}
The vertex and edge counts follow from~\eqref{eq:counts}. Let $H\subseteq G(r,t)$ be $C_6$-connected. If $e(H)\le1$, the desired bound follows from $t^2/Q\ge1$. Assume $e(H)\ge2$ and discard isolated vertices. By Lemma~\ref{lem:bipartite}, $\delta(H)\ge2$ and every two distinct vertices in the same bipartition class of $H$ have a common neighbor.

Define the support projections by $\pi_L(i,a)=[\lambda_i,a]$ and $\pi_R(j,b)=(\xi_j,b)$. Each map sends a copied vertex to the point or line that it represents. Let $J\subseteq B_Q$ consist of the images of the edges in $E(H)$ and their endpoints. Thus $J$ has no isolated vertices and at least one edge. It is pair-covered: if $z,z'$ are distinct vertices in the same class of $J$, choose preimages $u,u'\in V(H)$. These are distinct vertices in the same class of $H$, so they have a common neighbor $w$. The image of $w$ is then a common neighbor of $z,z'$ in $J$.

\smallskip\noindent\emph{Repeated supports: $t\ge Q$.}
Among the first $t$ terms of each periodic sequence, each element of $K$ occurs at most $\rho_0=\lceil t/Q\rceil\le2t/Q$ times. Lemma~\ref{lem:support} gives $e(J)\le3Q$. An edge that maps to the support edge $[\lambda,a]\sim(\xi,b)$ must have indices satisfying $\lambda_i=\lambda$ and $\xi_j=\xi$. The intercept $a$ and second coordinate $b$ are already fixed. Thus at most $\rho_0^2$ edges map to this support edge, and
\[
e(H)\le3Q\rho_0^2\le12\frac{t^2}{Q}.
\]

\smallskip\noindent\emph{Injective supports: $r\le t<Q$.}
Because $t<Q$, the first $t$ terms of each sequence lie within the first period and are distinct. Hence both support projections are injective on their entire vertex classes: distinct vertices have distinct images. It follows that $H\cong J$. In particular, $\delta(J)=\delta(H)\ge2$, so Lemma~\ref{lem:support} shows that $H$ is the incidence graph of a near-pencil. Let its exceptional point be $P_0=(x_0,y_0)$ and its long line be $\ell_0=[s_0,a_0]$.

Call each three-edge path from $P_0$ to $\ell_0$ a \emph{petal}. Consider a petal $P_0-\ell-P-\ell_0$, where $\ell=[s,a]$ and $P=(x,y)$. The line passes through $P_0$, so $a=y_0-sx_0$. The point $P$ lies on both $\ell$ and $\ell_0$, so $sx+a=s_0x+a_0$. Substituting for $a$ gives
\begin{equation}\label{eq:petal}
(s-s_0)(x-x_0)=s_0x_0+a_0-y_0=:c.
\end{equation}
The constant $c$ is nonzero because $P_0\notin\ell_0$. Moreover, distinct petals give distinct pairs $(s,x)$. Indeed, $s$ determines the unique line of that slope through $P_0$, and $x$ determines the unique point of first coordinate $x$ on $\ell_0$.

Every petal slope $s$ is one of $\lambda_1,\ldots,\lambda_t$, and the first coordinate $x$ of every petal point is one of $\xi_1,\ldots,\xi_t$. The first $t$ terms of either sequence lie in at most $D_0=\lceil t/r\rceil\le2t/r$ layers. By Lemma~\ref{lem:intersection}, applied to~\eqref{eq:petal}, and the fact that distinct petals give distinct pairs, the number $u$ of petals satisfies
\[
u\le4D_0^2\le16\frac{t^2}{r^2}=16\frac{t^2}{Q}.
\]
A near-pencil with $u$ petals has $3u$ edges, so $e(H)\le48t^2/Q$. The two cases cover every $t\ge r$.
\end{proof}

We choose $Q\asymp n^\beta$ and $t\asymp n^{1-\beta}$, giving $2tQ\asymp n$ vertices and $t^2Q\asymp n^{2-\beta}$ edges. The construction bound then has order $t^2/Q\asymp n^{2-3\beta}$. The next proof supplies the rounding details and explicit constants.

\begin{proposition}\label{prop:parameters}
For every fixed $0<\beta\le2/3$ and all sufficiently large $n$, there is an $n$-vertex graph with $\lfloor n^{2-\beta}\rfloor$ edges in which every $C_6$-connected subgraph has at most $580800n^{2-3\beta}$ edges.
\end{proposition}

\begin{proof}
Put $X=n^{\beta/2}$ and $\kappa=\lfloor X/9\rfloor$. For all sufficiently large $n$, we have $\kappa\ge2$, so Bertrand's postulate gives a prime $r$ with $\kappa<r<2\kappa$. Since $\kappa\to\infty$, we may also assume that $r$ is odd and greater than $4$. If $X\ge90$, then
\[
\kappa\ge\frac{X}{9}-1\ge\frac{X}{10},
\qquad \frac{X}{10}<r<\frac{2X}{9}.
\]
Put $Q=r^2$ and $t=\lceil10n^{1-\beta}\rceil$. Then
\[
\frac{n^\beta}{100}<Q<\frac{4n^\beta}{81},
\qquad
\frac{t}{r}>\frac{10n^{1-\beta}}{2n^{\beta/2}/9}
=45n^{1-3\beta/2}\ge45,
\]
so $t\ge r$. Moreover,
\[
2tQ<2(10n^{1-\beta}+1)\frac{4n^\beta}{81}
=\frac{80}{81}n+\frac{8}{81}n^\beta<n
\]
for large $n$, whereas
\[
t^2Q>(10n^{1-\beta})^2\frac{n^\beta}{100}=n^{2-\beta}.
\]
Delete edges until exactly $\lfloor n^{2-\beta}\rfloor$ remain, and add isolated vertices to obtain an $n$-vertex graph $G'$. Every $C_6$-connected subgraph $H\subseteq G'$ lies in $G(r,t)$ after isolated vertices are removed. Proposition~\ref{prop:construction} and $t\le11n^{1-\beta}$ therefore give
\[
e(H)\le48\frac{t^2}{Q}<580800n^{2-3\beta}.
\]
\end{proof}

\subsection{The lower bound and completion of the proof}\label{sec:lower}
The lower-bound argument in~\cite[Section~2]{DER1984} counts four-cycles and applies in the denser range. Here we count three-edge paths, which gives the estimate needed below $n^{3/2}$ edges.

\begin{proposition}\label{prop:lower}
Let $G$ be an $n$-vertex graph with $m$ edges. If $m\ge8n$ and $m^3/(64n^4)\ge2$, then $G$ contains a $C_6$-connected subgraph with at least $m^3/(64n^4)$ edges.
\end{proposition}

\begin{proof}
Choose a bipartite subgraph $G_\times=(A,B;E_\times)$ with $q=e(G_\times)\ge m/2$, by taking a partition with at least the expected number of crossing edges. Set $d_0=q/(2n)\ge2$ and repeatedly delete vertices of current degree less than $d_0$. The deletions remove fewer than $nd_0=q/2$ edges. Thus the remaining graph $J=(A',B';E')$ satisfies
\begin{equation}\label{eq:prune}
M_0:=e(J)>\frac q2,\qquad \delta(J)\ge\frac{q}{2n}\ge2.
\end{equation}
For $x\in A'$ and $y\in B'$, let $p(x,y)$ count the simple three-edge paths from $x$ to $y$. Counting by the middle edge gives
\begin{equation}\label{eq:paths}
\begin{aligned}
\sum_{x\in A'}\sum_{y\in B'}p(x,y)
&=\sum_{uv\in E(J)}(d_J(u)-1)(d_J(v)-1)\\
&\ge\frac14 M_0\delta(J)^2>\frac{q^3}{32n^2}.
\end{aligned}
\end{equation}
Here $d_J(z)-1\ge d_J(z)/2$ for every vertex $z$. For each middle edge $uv$, excluding $u$ and $v$ from the choices of endpoints ensures that all four vertices are distinct. Since $|A'||B'|\le n^2/4$, some pair $(x,y)$ satisfies
\begin{equation}\label{eq:average}
p(x,y)\ge\frac{q^3}{8n^4}\ge\frac{m^3}{64n^4}\ge2.
\end{equation}

Let $U=N_J(x)\setminus\{y\}$ and $V=N_J(y)\setminus\{x\}$, and let $J_0$ consist of the edges between $U$ and $V$, with isolated vertices removed. Continue to write $U,V$ for its remaining vertex classes. Its edges correspond exactly to the middle edges of the paths counted by $p(x,y)$, so $a:=e(J_0)=p(x,y)$. Adjoin $x,y$ and the edges from $x$ to $U$ and from $V$ to $y$; all these edges belong to $J$.

If $J_0$ is a star, its center and the opposite added vertex form a $K_{2,a}$ with its leaves. This graph is $C_4$-connected and has $2a$ edges. Otherwise, $|U|,|V|\ge2$, and the whole enlarged graph has minimum degree at least two. It is bipartite with classes $\{x\}\cup V$ and $\{y\}\cup U$. Two vertices of $U$ share $x$, and two vertices of $V$ share $y$. Also, $x$ and any $v\in V$ share a neighbor in $U$, while $y$ and any $u\in U$ share a neighbor in $V$. Lemma~\ref{lem:bipartite} therefore shows that this graph is $C_6$-connected. In either case we obtain at least $a\ge m^3/(64n^4)$ edges.
\end{proof}

\begin{proof}[Proof of Theorem~\ref{thm:sparse}]
Let $m=\lfloor n^{2-\beta}\rfloor$. If $0<\beta<2/3$, then $m/n\to\infty$, and
\[
\frac{m^3}{64n^4}=\left(\frac1{64}+o(1)\right)n^{2-3\beta}\longrightarrow\infty.
\]
Proposition~\ref{prop:lower} gives the lower bound, and Proposition~\ref{prop:parameters} gives the upper bound. At $\beta=2/3$, the latter proposition gives the upper bound $580800$, while a one-edge subgraph gives the lower bound $1$.
\end{proof}

\section{Proof of Theorem~\ref{thm:dense}}\label{sec:dense}

We use the following consequence of Duke, Erd\H{o}s, and R\"odl's characterization \cite[Proposition~0, p.~264]{DER1992}: deleting isolated vertices from a $C_4$-connected graph leaves a complete multipartite graph, and $K_{a,b}$ is $C_4$-connected whenever $a,b\geq2$. The first assertion also holds for a graph with one edge. We call a complete bipartite subgraph a \emph{biclique}, and a pair of disjoint vertex sets with no edges between them a \emph{bipartite hole}.

Write $G=K_n-R$, where $R$ has $s$ edges, and let $E_R(A,B)$ denote the edges of $R$ between disjoint sets $A,B$. Set
\begin{equation}\label{eq:dense-parameters}
d=\frac{2s}{n},\qquad M=\frac{n^2}{d}=\frac{n^3}{2s},\qquad
k=\left\lfloor\frac{n}{d}\right\rfloor.
\end{equation}
The assumptions give $d\to\infty$, $d=o(n)$, $k\to\infty$, and $n=o(M)$. All estimates below hold along the given sequence $s=s(n)$.

\paragraph{The lower bound.}
We use the star-counting argument in \cite[proof of Theorem~1, p.~264]{DER1992}. It gives a $k$-set $X$ with $b$ common neighbors in $G$ such that
\begin{equation}\label{eq:jensen}
b\geq\frac{\sum_{v\in V(G)}\binom{d_G(v)}{k}}{\binom{n}{k}}
\geq n\frac{\binom{n-1-\lceil d\rceil}{k}}{\binom{n}{k}}.
\end{equation}
The second inequality uses convexity of the binomial coefficients in the degree and rounds the average degree down. The fixed-$k$ asymptotic in the cited theorem does not apply directly here, so we give the needed uniform estimate. Put $d_+=\lceil d\rceil$. Since $k/n=O(1/d)$ and $(d_++1)/(n-k)=o(1)$, expansion of the logarithm gives
\begin{equation}\label{eq:log-product}
\begin{aligned}
\log\frac{\binom{n-1-d_+}{k}}{\binom{n}{k}}
&=\sum_{\ell=0}^{k-1}\log\left(1-\frac{d_++1}{n-\ell}\right)\\
&=-(d_++1)\sum_{\ell=0}^{k-1}\frac{1}{n-\ell}
+O\left((d_++1)^2\sum_{\ell=0}^{k-1}\frac{1}{(n-\ell)^2}\right)\\
&=-1+O\left(\frac1d+\frac dn\right).
\end{aligned}
\end{equation}
Indeed, the two sums are $k/n+O(k^2/n^2)$ and $O(k/n^2)$, respectively. Their error contributions are $O((d+1)k^2/n^2)=O(1/d)$ and $O((d+1)^2k/n^2)=O(d/n)$, while $(d_++1)k/n=1+O(1/d+d/n)$. Thus $b\geq(e^{-1}-o(1))n$. The set $X$ and its common neighbors are disjoint and form a biclique. Since $kn=M+O(n)=(1+o(1))M$,
\begin{equation}\label{eq:biclique-lower}
kb\geq(e^{-1}-o(1))M.
\end{equation}
Both parts have at least two vertices for large $n$, which proves the lower bound.

\paragraph{The upper bound.}
We follow the strategy of the upper-bound proof of \cite[Theorem~1, pp.~265--267]{DER1992}: first bound bicliques in a random graph, then use those bounds to exclude large complete multipartite subgraphs. The estimates below apply uniformly as the edge density tends to one.

Fix $\varepsilon>0$. Choose a fixed $0<\delta<1$ so small that
\begin{equation}\label{eq:delta}
(1-\delta)(1/e+\varepsilon)>1/e+\varepsilon/2,
\end{equation}
and set
\begin{equation}\label{eq:random-parameters}
\gamma=\frac{1/e+\varepsilon}{2},\qquad
D=2\sqrt{\frac{1/e+\varepsilon}{\delta}},\qquad
\rho=d^{-1/4},\qquad p=(1-\rho)\frac{d}{n-1}.
\end{equation}
Because $d\to\infty$ and $d=o(n)$, we have $0<p<1$ for all sufficiently large $n$. Let $R_0\sim G(n,p)$ and $G_0=K_n-R_0$. We call $|A|+|B|$ the support size of a biclique $K_{A,B}$; this is simply its total number of vertices. The next lemma gives one bound for all bicliques and a stronger bound for those with small support.

\begin{lemma}\label{lem:holes}
With probability tending to one, all three of the following assertions hold:
\begin{enumerate}
\item[(i)] every biclique in $G_0$ has fewer than $(1/e+\varepsilon/2)M$ edges;
\item[(ii)] every biclique in $G_0$ with support size at most $Dn/\sqrt{d}$ has fewer than $\gamma M$ edges;
\item[(iii)] $e(R_0)\leq s$.
\end{enumerate}
\end{lemma}

\begin{proof}
\emph{Assertion (i).} For a fixed $a$-set $A$, let $Z_A=\{v\notin A:E_{R_0}(\{v\},A)=\varnothing\}$. Then $|Z_A|\sim\operatorname{Bin}(n-a,(1-p)^a)$. Each vertex outside $A$ belongs to $Z_A$ with probability $(1-p)^a$, and these events are independent because they concern disjoint sets of possible edges.

Put $T=(1/e+\varepsilon/2)M$ and $L=\lceil\sqrt{T}\rceil=O(n/\sqrt{d})$. Any biclique $K_{U,W}$ with $|U||W|\geq T$ contains a biclique with one side $A$ of size $a\leq L$ and at least $T$ edges. Indeed, if neither side has size at most $L$, take $a=\lceil T/|W|\rceil$ vertices of $U$. Then $a\leq|U|$, and $T\leq L^2$ with $|W|>L$ gives $a\leq L$. Thus it is enough to rule out $a$-sets $A$, with $1\leq a\leq L$, for which $a|Z_A|\geq T$.

For $0\leq u,v\leq 1$, write
\[
D_{\mathrm{KL}}(u\mathbin\Vert v)
=u\log\frac{u}{v}+(1-u)\log\frac{1-u}{1-v},
\]
with the usual continuous conventions at the boundary. Fix such an $A$ and put $x=pa$, $q_a=(1-p)^a$, and $\theta_a=T/[a(n-a)]$. Choose a constant $C$ with $1/e<C<1/e+\varepsilon/2$. Uniformly for $a\leq L$, we have $a=o(n)$ and
\[
\frac{pT}{n-a}
=\left(\frac{1}{e}+\frac{\varepsilon}{2}\right)(1-\rho)
\frac{n^2}{(n-1)(n-a)}
=\frac{1}{e}+\frac{\varepsilon}{2}+o(1).
\]
Consequently, for all sufficiently large $n$,
\begin{equation}\label{eq:theta}
\theta_a\geq\frac{C}{pa}=\frac{C}{x}.
\end{equation}
If the event $a|Z_A|\geq T$ can occur, then $\theta_a\leq 1$ and hence $x\geq C$. Furthermore,
\[
q_a\leq e^{-x}<\frac{C}{x}\leq\theta_a,
\]
because $xe^{-x}\leq1/e<C$. For $u>v$, the function $D_{\mathrm{KL}}(u\mathbin\Vert v)$ increases with $u$ and decreases with $v$. The binomial Chernoff bound and this monotonicity give
\[
\begin{aligned}
\Pr(a|Z_A|\geq T)
&\leq\exp\bigl(-(n-a)D_{\mathrm{KL}}(\theta_a\mathbin\Vert q_a)\bigr)\\
&\leq\exp\bigl(-(n-a)D_{\mathrm{KL}}(C/x\mathbin\Vert e^{-x})\bigr).
\end{aligned}
\]
Set $I_C(x)=D_{\mathrm{KL}}(C/x\mathbin\Vert e^{-x})$ for $x\geq C$. The inequality $xe^{-x}<C$ shows that $I_C(x)>0$. Moreover, $I_C(C)=C$ and
\[
I_C(x)
=\frac{C}{x}\left(x+\log\frac{C}{x}\right)
+\left(1-\frac{C}{x}\right)\log\frac{1-C/x}{1-e^{-x}}
=C+o(1)
\]
as $x\to\infty$. Continuity and positivity on each compact interval, together with this positive limit, give $\kappa_0:=\inf_{x\geq C}I_C(x)>0$. Since $a\leq L=o(n)$,
\[
\Pr(a|Z_A|\geq T)\leq\exp(-\kappa_0 n/2)
\]
for all sufficiently large $n$. Since $L=o(n)$, we have $L\leq n/2$ for large $n$, and the number of possible sets is at most
\[
\sum_{a=1}^{L}\binom{n}{a}
\leq(L+1)\binom{n}{L}
\leq(L+1)\left(\frac{en}{L}\right)^L
=\exp\left(O\left(\frac{n\log d}{\sqrt{d}}\right)\right)
=e^{o(n)}.
\]
A union bound proves (i), since the total probability tends to zero.

\emph{Assertion (ii).} Put $S=Dn/\sqrt{d}=o(n)$ and $S'=\lceil S\rceil$. To count ordered disjoint pairs $(A,B)$ with $|A|+|B|\leq S$, choose their union and assign each vertex to one of the two sets. Since $S'=o(n)$, the summands below increase throughout the range, giving at most
\[
\begin{aligned}
\sum_{u\leq S'}\binom{n}{u}2^u
&\leq(S'+1)\binom{n}{S'}2^{S'}\\
&\leq(S'+1)\left(\frac{2en}{S'}\right)^{S'}
=\exp\left(O\left(\frac{n\log d}{\sqrt{d}}\right)\right)
=e^{o(n)}.
\end{aligned}
\]
For a fixed pair with $|A||B|\geq\gamma M$,
\[
\begin{aligned}
\Pr(E_{R_0}(A,B)=\varnothing)
&=(1-p)^{|A||B|}\\
&\leq\exp(-p\gamma M)\\
&=\exp\left(-\gamma(1-\rho)\frac{n^2}{n-1}\right)
=\exp\bigl(-(\gamma-o(1))n\bigr).
\end{aligned}
\]
Another union bound proves (ii).

\emph{Assertion (iii).} Let $Z=e(R_0)$ and $\mu=\mathbb{E}Z=(1-\rho)s$. Writing $s=(1+\alpha)\mu$, we have $\alpha=\rho/(1-\rho)<1$ for large $n$. The multiplicative Chernoff bound yields
\[
\Pr(Z>s)\leq\exp\left(-\frac{\alpha^2\mu}{3}\right)
\leq\exp\left(-\frac{\rho^2s}{3}\right).
\]
Since $\rho^2s=d^{-1/2}s=n\sqrt{d}/2\to\infty$, the right-hand side tends to zero, proving (iii).
\end{proof}

\begin{proof}[Completion of the upper bound]
Choose an outcome of $R_0$ satisfying Lemma~\ref{lem:holes}. Since $e(R_0)\leq s<\binom{n}{2}$ for all sufficiently large $n$, we can add edges to obtain a graph $R$ with exactly $s$ edges. This only deletes edges from $G_0$, so both biclique bounds still hold in $G=K_n-R$.

Suppose that $G$ contains a complete multipartite subgraph with at least $T_0=(1/e+\varepsilon)M$ edges. Choose such a subgraph $H$ with as few vertices as possible, and put $E=e(H)$. Deleting any vertex leaves fewer than $T_0$ edges. Since every degree is at most $n-1$,
\[
T_0\leq E<T_0+n-1.
\]
Write the part sizes in decreasing order as $t_1=\tau\geq t_2\geq\cdots\geq t_q$. Put $r_0=\sum_{i\geq 2}t_i$ for the number of vertices outside the largest part and $\sigma=\tau+r_0$ for the total number of vertices. Let $E_{\mathrm{rest}}=\sum_{2\leq i<j\leq q}t_it_j$ count the edges between parts other than the largest. Then $E=\tau r_0+E_{\mathrm{rest}}$.

If $E_{\mathrm{rest}}\leq\delta E$, the largest part and the union of all other parts form a biclique with
\[
\tau r_0=E-E_{\mathrm{rest}}\geq(1-\delta)E
\geq(1-\delta)T_0>(1/e+\varepsilon/2)M
\]
edges, contradicting Lemma~\ref{lem:holes}(i). Otherwise $E_{\mathrm{rest}}>\delta E$. Since $E_{\mathrm{rest}}\leq r_0^2/2$, we obtain
\begin{equation}\label{eq:r0}
r_0>\sqrt{2\delta E}.
\end{equation}
As every $t_i\leq\tau$, we have $\sum_i t_i^2\leq\tau\sigma$. Consequently,
\[
E=\frac{\sigma^2-\sum_i t_i^2}{2}
\geq\frac{\sigma^2-\tau\sigma}{2}
=\frac{\sigma r_0}{2}.
\]
Combining the last two displays gives
\[
\sigma\leq\frac{2E}{r_0}<\sqrt{\frac{2E}{\delta}}.
\]
Since $n=o(M)$ and $E<T_0+n$,
\[
E\leq\left(\frac{1}{e}+\varepsilon+o(1)\right)\frac{n^2}{d}.
\]
Consequently,
\[
\sqrt{\frac{2E}{\delta}}
\leq\sqrt{\frac{2(1/e+\varepsilon+o(1))}{\delta}}\frac{n}{\sqrt{d}}
<D\frac{n}{\sqrt{d}}
\]
for all sufficiently large $n$, by the definition of $D$ in~\eqref{eq:random-parameters}. Hence
\begin{equation}\label{eq:support}
\sigma<D\frac{n}{\sqrt{d}}.
\end{equation}
Color each part of $H$ independently red or blue, choosing each color with probability $1/2$. Let $A$ and $B$ be the unions of the red and blue parts. Each edge joins parts of different colors with probability $1/2$, so some coloring has at least $E/2$ edges between $A$ and $B$. Since every part has one color, $A$ and $B$ are unions of different parts and form a biclique in $G$. Its support size is at most $\sigma<Dn/\sqrt{d}$, and it has at least $E/2\geq(1/e+\varepsilon)M/2=\gamma M$ edges. This contradicts Lemma~\ref{lem:holes}(ii).

Thus $G$ contains no complete multipartite subgraph with at least $T_0$ edges. The structural result cited at the start of this section gives the same conclusion for $C_4$-connected subgraphs. Since $\varepsilon>0$ is arbitrary, this upper bound and~\eqref{eq:biclique-lower} give $f_2\bigl(n,\binom{n}{2}-s\bigr)=(1/e+o(1))M$. Substituting $M=n^3/(2s)$ completes the proof.
\end{proof}

\section*{Acknowledgments}
This research was supported by the National Key R\&D Program of China (Grant No.\ 2023YFA1010202), the National Natural Science Foundation of China (Grant No.\ 12401448), and the Natural Science Foundation of Fujian Province (Grant No.\ 2024J08030).

\section*{Declaration on the use of generative AI}

The authors used generative AI tools to assist in discussing proof strategies, checking proofs, and
improving the exposition.
This project has been under investigation since 2022. Some of the key ideas were contributed by the authors. 
The authors take full responsibility for the mathematical arguments, results,
and conclusions, all of which were carefully reviewed and verified by them.

\end{document}